\documentclass{notices}
\usepackage{amsfonts,amssymb,amsmath,amscd}

\usepackage{graphicx}
\usepackage{hyperref, url}
\usepackage{color}
\usepackage{dsfont}

\newtheorem{thm}{Theorem} 
\newtheorem{prop}[thm]{Proposition}

\newtheorem{lemma}[thm]{Lemma}
\newtheorem{defn}[thm]{Definition}

\newcommand{\dx} {\; \mathrm{d} x}
\newcommand{\dd} {\; \mathrm{d}}
\newcommand{\ddsw} {\; \mathrm{d}\sigma\mathrm{d} w}
\newcommand{\ddvw} {\; \mathrm{d}v\mathrm{d}w}
\newcommand{\ddsrz} {\; \mathrm{d}\sigma\mathrm{d}r\mathrm{d}z}

\DeclareMathOperator{\dv}{div}

\title{The Landau equation and Fisher information}
\date{}

\author{
  Nestor Guillen
  \affil{The first author is a professor of mathematics at Texas State University. His email address is nestor@txstate.edu.}
  \and
  Luis Silvestre
  \affil{The second author is a professor of mathematics at The University of Chicago. His email address is luis@math.uchicago.edu.}
}

\begin{document}

\maketitle

\begin{figure}[t]

\centering
\includegraphics[scale=0.25]{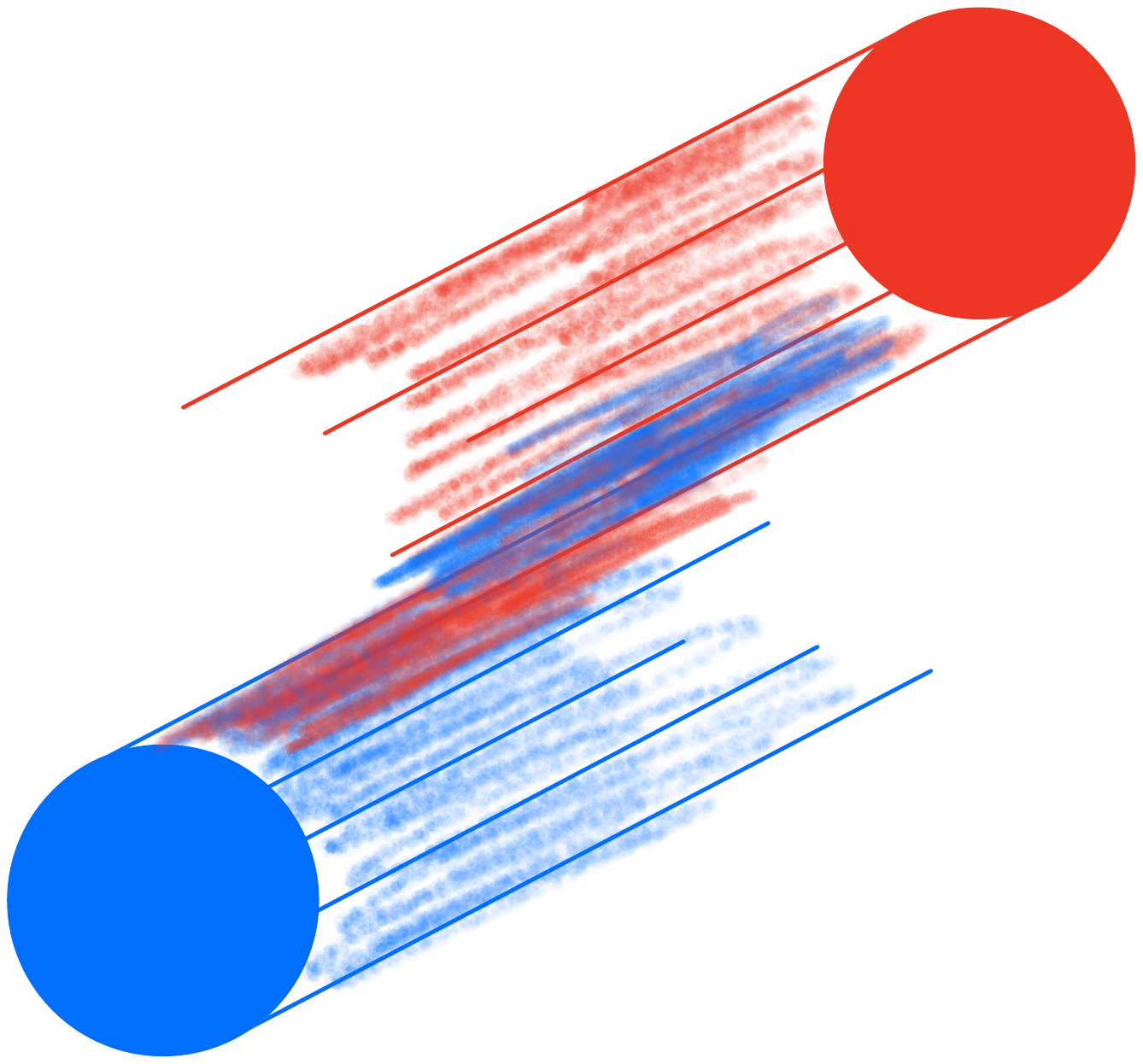}

\vspace{-0.2in}

\end{figure}

\section{Introduction}

The following initial value problem is of great interest in kinetic theory: to find, given an initial mass density $f_{\text{in}}:\mathbb{R}^3\to\mathbb{R}$,  a function $f(t,v)$ smooth for all times $t>0$ and velocities $v\in\mathbb{R}^3$, which starts at time $t=0$ from $f_{\text{in}}$ and solves for $t>0$
\begin{align}\label{e:Landau Coulomb}
  \partial_t f = \bar a_{ij} \partial_{ij}^2f + f^2	 \text{ where } \bar a := -D^2(-\Delta)^{-2}f
\end{align}	
This equation is known as the homogeneous Landau equation with Coulomb potential. 
It is an equation with a rich structure, it has both integro-differential and nonlinear terms. The term $f^2$ raises the possibility that solutions starting from a smooth initial $f_{\text{in}}$ might develop singularities (``blow up'') in final time.

In this note we explain how the blow up question for the Landau equation was ruled out through the study of the Fisher information, which turns out to be monotone in time. In \cite{GuiSil2023} we show that if $f$ is a classical solution of \eqref{e:Landau Coulomb} then
\begin{align*}
  \frac{d}{dt}\int_{\mathbb{R}^3}\frac{|\nabla f|^2}{f}\dd v\leq 0.
\end{align*}
This inequality is obtained via a new ``lifted equation'' taking place in $\mathbb{R}^6$. This lifted equation amounts to the heat equation on the sphere, a connection that was at first somewhat surprising but quite natural in hindsight. As this equation is local and linear, it is much easier to analyze than \eqref{e:Landau Coulomb}, and yet it encodes enough about \eqref{e:Landau Coulomb} to reach conclusions about how the Fisher information evolves with time. 

\subsection{Origin of the Landau equation} It is proper to start by reviewing the Boltzmann equation and how the Landau equation \eqref{e:Landau Coulomb} arises from it. A  complete introduction to the Boltzmann equation and its origins can be found in Cercignani's book \cite{Cercignani1975book}. A thorough survey of the mathematical aspects of the theory and the state of the field at the turn of the century is in Villani's review \cite{villani2002review}. 

In 1872 Boltzmann \cite{Boltzmann1872} derived an equation to model the evolution of the statistics of particles in a monoatomic gas, it takes the form
\begin{align}\label{e:Boltzmann}
\partial_tf + v \cdot \nabla_x f = q(f)
\end{align}
The function $f = f(t,x,v)$ represents the density of particles with position $x$ and velocity $v$, the equation indicates how $f$ changes with time as a result of the inertia of particles (this is captured by the term $v\cdot \nabla_xf$) and the result of interactions between particles, such as collisions (captured by the term $q(f)$). The \emph{collision operator} $q(f)$ introduced by Boltzmann should be thought of as providing a mass balance for how many particles are \emph{gaining} the velocity $v$ minus how many are \emph{losing} the velocity $v$ to end at another velocity as a result of a collision.

Boltzmann showed that the collision term $q(f)$ evaluated $(t,x,v)$ is given by the integral (we write $f(v)$ for $f(t,x,v)$, and $f(w)$ for $f(t,x,w)$ et cetera)

\vspace{-0.2in}
\begin{align}\label{e:Boltzmann collision operator}
\int_{\mathbb{R}^3}\int_{\mathbb{S}^2} (f(v_\sigma)f(w_\sigma)-f(v)f(w))B\ddsw
\end{align}
Here, 
\vspace{-0.1in}
\begin{align*}
  B & =  \alpha(|v-w|)b( \sigma \cdot \tfrac{v-w}{|v-w|})\\
  v_\sigma & = \tfrac{1}{2}(v+w)+\tfrac{1}{2}|v-w|\sigma\\
  w_\sigma & = \tfrac{1}{2}(v+w)-\tfrac{1}{2}|v-w|\sigma
\end{align*}
Here we will refer to $\alpha:(0,\infty)\to [0,\infty)$ as the \emph{potential}. Maxwell \cite{Maxwell1867} computed $\alpha(r)$ and $b(\cos(\theta))$ in terms of the force modeling interactions (collisions) between pairs of particles. For a central force given by $s$-th power of the inverse distance we have
\begin{align*}
  \alpha(r) = r^{\gamma},\textnormal{ where } \gamma := \frac{s-5}{s-1}
\end{align*}
For the case of Coulomb interactions this means $s=2$ and so $\gamma=-3$. It is well known that
\begin{align*}
  b(\sigma\cdot \sigma_0) \approx |\sigma-\sigma_0|^{-2-\frac{2}{s-1}}
\end{align*}
Here $\sigma_0 = \tfrac{v-w}{|v-w|}$. The integral in \eqref{e:Boltzmann collision operator} is not finite for any non-constant $f$ when $s>2$. One way of thinking about this is that when $s$ is close to $2$ then $b(\sigma\cdot\sigma_0)$ assigns most of its weight to those collisions with $\sigma \approx \sigma_0$, what are known as \emph{grazing collisions}.

For this reason it is interesting to consider the asymptotic regime where the contributions to the collisional integral are concentrated on directions $\sigma$ close to $\sigma_0 = \tfrac{v-w}{|v-w|}$. This is known in the literature as the \emph{grazing limit}. In this regime, most of the particle interactions produce asymptotically small changes in velocity. It is a regime of physical relevance because it is used to describe particles that repel each other by Coulombic forces. 

The Landau collision operator is obtained from the Boltzmann collision operator in the grazing limit. Different families of collision kernels can be considered, with the common feature that they concentrate near $\sigma = \sigma_0$. The inner spherical integral in \eqref{e:Boltzmann collision operator} is an spherical integro-differential operator applied to the function $f(v)f(w)$ with the kernel $b(\sigma\cdot \sigma_0)$, which is a function of the angle between $\sigma$ and $\sigma_0$. As $b$ concentrates on small angles, the integro-differential operator converges to the Laplace-Beltrami operator on the sphere $\mathbb{S}^2$ with respect to the variable $\sigma$, applied to the function $f(v)f(w)$,
\begin{equation} \label{e:funny-Landau}
  q_L(f)(v) = \int_{\mathbb{R}^3} \alpha(|v-w|) \Delta_\sigma[f(v)f(w)] \dd w.
\end{equation} 
It is arguably not obvious that when $\alpha(r)=r^{-3}$, this operator is the same Landau collision operator that we wrote at the beginning of this article \eqref{e:Landau Coulomb}. The Landau collision operator has several equivalent formulations. It is useful to know its many faces, since they can be used for different purposes.

For general $\alpha(r)$, the following are  common representations of the Landau collision operator 
\begin{itemize}
\item{\bf Nondivergence form}
\[ q_L(f) = \bar a_{ij} \partial_{ij} f + \bar c f,\]
where $\bar c = -\partial_{ij} \bar a_{ij}$ and $\bar a_{ij} = a_{ij} \ast f$, for the matrix-valued function $a_{ij}$ given by
\begin{equation} \label{e:aij}
a_{ij}(z) = \alpha(|z|) \left(|z|^2 \delta_{ij} - z_i z_j \right)
\end{equation}
\item{\bf Divergence form}
\[ q_L(f) = \partial_i \left( \bar a_{ij} \partial_j f - \bar b_i f \right),\]
where $\bar a_{ij}$ is as before and $\bar b_i = \partial_j \bar a_{ij}$.
\item{\bf Integral form}
\begin{align*} &q_L(f) = \\
  &\partial_i \int_{\mathbb{R}^3} a_{ij}(v-w) \left( f(w) \partial_j f(v) - f(v) \partial_j f(w) \right) \mathrm{d} w.
\end{align*}
\end{itemize}

It is not too complicated to verify that the expression in \eqref{e:funny-Landau} coincides with the integral form of the Landau collision operator. The formulas for the operator in divergence and nondivergence form are computed from this integral formulation in an attempt to make the Landau collision operator look like a classical diffusion operator.

\subsection{The homogeneous regime and the blow up question}

The study of kinetic equations faces two main analytical difficulties. First, while the collision operator has a diffusive regularization effect, it is only with respect to the $v$-variable. To get a regularization in all variables we have to understand the interaction between the diffusion in $v$ and the transport in $x$. This is where hypoelliptic theory applies to kinetic equations. The other difficulty is that the collision operator $q(f)$ has a rather complex structure --it is a nonlocal and quadratic expression in $f$.

For configurations that are homogeneous in space, namely when $f(t,x,v)=f(t,v)$, the transport effects in the equation vanish, provides us with an ideal situation to analyze the effects of collisions. In this ``spatially homogeneous'' regime the Boltzmann equation \eqref{e:Boltzmann} simplifies to
\begin{equation}\label{e:homogeneous Landau}
    \partial_tf = q(f)
\end{equation} 

The spatially homogeneous regime, as simplified as it might be, is still of considerable interest to physicists -- see for instance the discussion in Pauli's \emph{Statistical Mechanics}  \cites{Pau1973vol4}[Chapter 1]. Mathematically, it is a necessary first step before any attempt to understand the full equation \eqref{e:Boltzmann}.

The most natural idea to prove that the solution of the homogeneous Landau equation \eqref{e:homogeneous Landau} stays bounded and smooth for all time is to take advantage of the regularization effect of the diffusion term in $q_L(f)$ (either in divergence or non-divergence form). However, the presence of the nonlinear term $f^2$ in \eqref{e:Landau Coulomb} makes this approach difficult. This difficulty exists whenever the potential $\alpha(r)$ is too singular near $r=0$. The equation \eqref{e:Landau Coulomb} is reminiscent of the nonlinear heat equation $u_t = \Delta u + u^2$, which actually blows up in finite time. The question is then, what makes the Landau equation fundamentally different from the nonlinear heat equation?

The homogeneous problem \eqref{e:homogeneous Landau}, for either the Boltzmann or Landau collision operator, has a few well-known conserved quantities: the mass, momentum and energy. Moreover, the entropy of $f$ is monotone decreasing in time.
\begin{itemize}
\item {\bf Mass}
\[ M(f) = \int_{\mathbb{R}^3} f(v) \dd v\]
\item {\bf Momentum}
\[ P(f) = \int_{\mathbb{R}^3} v f(v) \dd v\]
\item {\bf Energy}
\[ E(f) = \int_{\mathbb{R}^3} \frac{|v|^2}{2} f(v) \dd v\]
\item {\bf Entropy}
\[ H(f) = \int_{\mathbb{R}^3} f(v) \ln f(v) \dd v \]
\end{itemize}

The conservation of the first three of these quantities, and the monotonicity of the entropy are classical results, obtained many years ago by Boltzmann himself \cite{Boltzmann1872}. These quantities, together with regularity estimates for parabolic equations, suffice to show that the solution to the space-homogeneous problem \eqref{e:homogeneous Landau} stays smooth for all time only when the potential $\alpha$ is not too singular. In the case of the Landau equation and $\alpha(r) = r^\gamma$, this is when $\gamma \geq -2$. The case of Coulombic potentials, which corresponds to $\gamma = -3$, is well outside this range. To rule out the possibility of a finite time blow up, we found a new monotone-in-time quantity, in addition to the entropy. This new quantity is the Fisher information.

The Fisher information is an integral functional over probability densities $f=f(v)$ and defined by
\begin{align*}
  I(f) = \int_{\mathbb{R}^N}\frac{|\nabla f|^2}{f}\;dv.
\end{align*}
McKean \cite{Mck1966} was the first to investigate the behavior of the Fisher information in the context of kinetic equations. Concretely, he showed that $i(f)$ is non-increasing in time for solutions of Kac's 1D model for the Boltzmann equation. McKean was skeptical this could be extended to the true Boltzmann equation. Several years later, Toscani \cite{Tos1992} showed the Fisher information was monotone decreasing for the Boltzmann equation in the Maxwell molecules case ($\alpha \equiv 1$). Villani extended Toscani's result to arbitrary dimension \cite{Vil1998}, still for the Boltzmann equation with Maxwell molecules. In that paper, Villani explains why their methods were restricted to the Maxwell molecules case $\alpha \equiv 1$, perhaps raising some skepticism about the possibility of extending the result to the case of Coulombic potentials. However, in \cite{villani2002review}*{Chapter 4, Section 3}, Villani reports that numerical tests performed by Buet and Cordier\cite{BuetCordier1998}[Section 5] suggest that the Fisher information may be monotone decreasing for a general potential $\alpha(r)$, including the Coulomb case. A quarter of a century later, and more than 50 years after McKean's original work, we were able to prove this is the case for the Landau equation in \cite{GuiSil2023}, and thus settle the question of global existence of global smooth solutions in that case. Later on,  Imbert, Villani, and the second author \cite{ImbSilVil2024} proved that the Fisher information is monotone decreasing for the Boltzmann equation for a large class of collision kernels that includes all the physically relevant cases of power law interactions. In this article, we review the ideas leading to the proof of the monotonicity of the Fisher information for the Landau equation.

\section{The structure of $q(f)$}

We discussed that the operator $q(f)$ can be written in different forms. When we write it in divergence or non-divergence form, it takes a particularly short form. Its integral formulation reveals some further structure. We will further reinterpret this integral formulation as a combination of the operations of tensor product of densities, a linear (elliptic) operator, and of projection onto a marginal.

We explain each of these terms. First, given functions $f,g:\mathbb{R}^3\to \mathbb{R}$ their tensor product $f\otimes g$ is the function in $\mathbb{R}^6 = \mathbb{R}^3\times\mathbb{R}^3$ given by
\begin{align}\label{e:definition tensor product}
  (f\otimes g) (v,w) := f(v)g(w)
\end{align}
Of particular relevance is the case when $g=f$, where $f\mapsto f\otimes f$ defines a map from mass densities in $\mathbb{R}^3$ to mass densities in $\mathbb{R}^6$. 

In the opposite direction we have the operation of computing a marginal density (thinking of $\mathbb{R}^6$ as $\mathbb{R}^3\times \mathbb{R}^3$), if $F:\mathbb{R}^3\times \mathbb{R}^3\to \mathbb{R}$, we define its projection $\pi F:\mathbb{R}^3\to \mathbb{R}$ by
\begin{align*}
\pi F(v) = \int_{\mathbb{R}^3}F(v,w)\dd w
\end{align*}
Then the Boltzmann collision operator has the form
\begin{align}\label{e:collision operator decomposition}
  q(f) = \pi (Q(f\otimes f))
\end{align}
where for (a sufficiently smooth) $F:\mathbb{R}^6\to\mathbb{R}$
\begin{align}\label{e:collision operator Boltzmann lifting}
Q(F)(v,w) := \int_{\mathbb{S}^2}(F(v_\sigma,w_\sigma)-F(v,w))B\dd\sigma
\end{align}
It is clear that $Q$ is a linear operator, in fact, it is an integro-differential (degenerate) elliptic operator.

Naturally, the same decomposition holds for the Landau operator. To find the appropriate operator $Q$ we analyze the integral form of the operator $q(f)$, which is given by
\begin{align*}
  &q(f)(v) = \\
  & \partial_i \int_{\mathbb{R}^3} a_{ij}(v-w) \left( f(w) \partial_j f(v) - f(v) \partial_j f(w) \right) \dd w \\
  &= \pi \left[ \partial_{v_i} a_{ij}(v-w) \left( f(w) \partial_j f(v) - f(v) \partial_j f(w) \right) \right] \\
  &= \pi \left[ (\partial_{v_i}-\partial_{w_i}) a_{ij}(v-w) (\partial_{v_j}-\partial_{w_j}) [f \otimes f]  \right]
\end{align*}
Thus, for the Landau operator we define 
\begin{align*}
Q(F) := (\partial_{v_i}-\partial_{w_i}) a_{ij}(v-w) (\partial_{v_j}-\partial_{w_j}) F
\end{align*}
and we have
\begin{align*}
  q(f) = \pi (Q(f\otimes f))
\end{align*}

Note that in the last line of the computation above we introduced, seemingly arbitrarily, the differentiation with respect to $w_i$. We are able to do it because the derivative with respect to $w_i$ of any function with appropriate decay integrates to zero with respect to $w$. Our purpose is to end up with an operator $Q$ that is elliptic, and that preserves symmetry: if $F(v,w) = F(w,v)$, then $Q(F)(v,w) = Q(F)(w,v)$.

\subsection{Lifting functionals}

The decomposition of $q(f)$ we have just described is meaningful because it reduces certain questions about the nonlinear Landau equation to questions about the linear diffusion operator $Q$.

Given $f(v)$, let us consider the following initial value problem in $(0,\infty)\times \mathbb{R}^6$,
\begin{align}\label{e:lifting equation Cauchy problem}
  \left \{ \begin{array}{rl} \partial_t F = Q(F) & \text{ for } t>0 \\ F = f\otimes f & \text{ at } t=0\end{array}\right.
\end{align}
We call this the \emph{lifted equation} and $Q(F)$ the \emph{lifted collision operator}. For the Landau equation, the lifted equation is a linear degenerate parabolic equation. For the Boltzmann equation, the lifted equation is a linear integro-differential diffusion.

If $f(t,v)$ is a smooth solution of the Landau equation and $F(t,v,w)$ is a smooth solution of \eqref{e:lifting equation Cauchy problem} with $F(0,v,w) = f(0,v)f(0,w)$ then
\begin{align*}
  (\partial_t f )_{\mid t=0} = (\partial_{t}\pi(F))_{\mid t=0}
\end{align*}
Put differently, the $v$-marginal to $F(t,\cdot)$, $t \mapsto \pi F(t,\cdot)$, describes a curve in the space of mass densities with first order contact with $t \mapsto f(t,\cdot)$ at time $t=0$. 

Thanks to this we are able to relate the derivative of a functional $j(f)$ along the flow of the (nonlinear) Landau equation with the derivative of a \emph{lifted} functional $J(F)$ along the flow of the lifted equation. This can be done for some special type of functionals that include the entropy and the Fisher information.

\begin{defn}\label{d:lifting}
 Consider two functionals, $j(f)$ and $J(F)$, defined respectively over probability densities in $\mathbb{R}^d$ and $\mathbb{R}^{2d}$. The functional $J(F)$ is said to be a ``lifting'' of the functional $j(f)$ if 
 
 1) If $F(v,w)=F(w,v)$ then $J(F) \geq 2j(\pi F)$.
  
 2) For any $f$, $J(f\otimes f) = 2j(f)$.
\end{defn}

\begin{prop} If $J$ is a lifting of $j$ then 
\begin{align*}
  2\langle j'(f),q(f)\rangle = \langle J'(f\otimes f),Q(f\otimes f)\rangle
\end{align*}

\end{prop}
Let us show at least ``one side'' of this equality, i.e. the inequality that is needed for the rest of our proof. Given a solution $f(t,v)$ starting from $f(v)$, consider the associated $F(t,v,w)$ solving the lifted equation and starting from $f\otimes f$. For every $t>0$ we have
\begin{align*}
  2j(\pi F(t)) \leq J(F(t))
\end{align*}
Given this is an equality for $t=0$ it follows that
\begin{align*}
  2 \frac{d}{dt}_{\mid t= 0} j(\pi F(t)) \leq \frac{d}{dt}_{\mit t=0}J(F(t))
\end{align*}
and the following inequality follows,
\begin{align*}
  2\langle j'(f),q(f)\rangle \leq \langle J'(f\otimes f),Q(f\otimes f)\rangle
\end{align*}
To prove the inequality in the other direction, we would have to flow the equation to negative time, which may be doable if $f$ is sufficiently smooth.

Finding a lifting $J$ for a functional $j$ is of great interest if we want to understand $\frac{d}{dt}j(f)$. Having such a $J$ furnishes us with a significant reduction of the problem of estimating the derivative: for the price of going up to twice the dimension, we only need to understand how $J(F)$ evolves in time for solutions of a linear degenerate parabolic equation, $\partial_tF =Q(F)$.


We are not aware of many other functionals besides $i(f)$ and $h(f)$ that have a natural lifting to $\mathbb{R}^{2d}$. Finding such functionals merits further investigation, since they would be amenable of the same type of analysis described in this article. A similar concept applies to distances between probability densities. The squared Wasserstein distance $W^2(f,g)$ is a functional that can be lifted to the space of probability densities in $\mathbb{R}^6$ since it satisfies the inequality
\[ W^2(F,G) \geq 2W^2(\pi F, \pi G),\]
with equality when $F=f \otimes f$ and $G=g\otimes g$.

\section{The Fisher information}\label{s:the Fisher information}

When we talk about the Fisher information functional for a density $f(v)$ we mean the integral
\begin{align}\label{e:Fisher}
  i(f) = \int_{\mathbb{R}^3}\frac{|\nabla f(v)|^2}{f(v)}\dd v
\end{align}
This functional has received considerable attention from mathematics since Linnik \cite{Linnik1959} used it as a key tool in his ``information theoretic'' proof of the Central Limit Theorem. The Fisher information, strictly speaking, is a more general notion in the original form introduced by Fisher \cite{Fisher1922} as part of his foundational work in statistics. In this original form, the Fisher information is associated to a parametric family of probability distributions $(f_\theta(x))_{\theta}$ in $\mathbb{R}^N$ and is given by
\begin{align*}
i(f,\theta)  = \int_{\mathbb{R}^N} \frac{(\partial_\theta f_\theta(v))^2}{f_\theta(v)}\dd v
\end{align*}
which is close in form to \eqref{e:Fisher}. An important class of examples are those where $f_\theta(v) = f(v+\theta e)$ for some $f$ and $e \in S^{d-1}$. Adding over the vectors $e$ of the canonical basis leads to $i(f)$ as in \eqref{e:Fisher}. More generally, if $b$ is a smooth, divergence-free vector field in $\mathbb{R}^N$ and $f_\theta$ denotes the push-forward $\Phi^{\theta}_{\#}f_0$ of a reference density $f_0$ with respect to the flow generated by $b$, then  one can check that
\begin{align*} 
  i(f,\theta=0) & = \int_{\mathbb{R}^N} \frac{(b(v)\cdot\nabla f_0(v))^2}{f_0(v)}\dd v
\end{align*}
For what follows, it is illustrative to see how \eqref{e:Fisher} evolves under the flow of the heat equation. If $f\geq 0$ is smooth and solves $\partial_tf=\Delta f$ it is well known that
\begin{align}\label{e:Fisher derivative heat equation}
  \frac{d}{dt}\int \frac{|\nabla f|^2}{f}\dd v = -2\int |D^2\ln f|^2f\dd v
\end{align}
This formula can be derived via integration by parts -although the derivation can feel magical. It is also well known that this is a computation that can be approached via Bochner's identity:
\begin{align}\label{e:Bochner Euclidean space}
 \Delta |\nabla g|^2 = 2|D^2g|^2+2\nabla g\cdot \nabla \Delta g
\end{align}
which holds for any smooth $g$ in Euclidean space. To see how \eqref{e:Fisher derivative heat equation} follows from \eqref{e:Bochner Euclidean space} we consider that
\begin{align*}
\frac{d}{dt}\int \frac{|\nabla f|^2}{f}\dx & = \int \partial_t \left ( \frac{|\nabla f|^2}{f} \right )\dd v
\end{align*}
together with the chain rule
\begin{align*}
  \partial_t \left ( \frac{|\nabla f|^2}{f} \right ) = 2\nabla \ln f\cdot \nabla (\Delta f)-2|\nabla \ln f|^2\Delta f
\end{align*}
Then we note the first term on the right hand side is almost the last one in \eqref{e:Bochner Euclidean space} with $g= \ln f$. In fact,
\begin{align*}
  \nabla \ln f\cdot \nabla \Delta f = f\nabla \ln f\cdot \nabla \Delta\ln f + \dv(|\nabla \ln f|^2\nabla f)
\end{align*}
This combined with the divergence theorem yields
\begin{align*}
\int f\nabla \ln f\cdot \nabla \Delta f\dd v=   \int f\nabla \ln f\cdot \nabla \Delta\ln f \dd v
\end{align*}
By Bochner's formula \eqref{e:Bochner Euclidean space}, the last integral equals
\begin{align*}
-2\int f|D^2\ln f|^2f\dd v + \int f\Delta (|\nabla \ln f|^2)\dd v
\end{align*}
Then,
\begin{align*}
\frac{d}{dt}\int \frac{|\nabla f|^2}{f}\dd v & = -2\int f|D^2\ln f|^2f\dd v\\
  & \;\;\;\;+\int f\Delta (|\nabla \ln f|^2)\dd v\\
  & \;\;\;\;-\int |\nabla \ln f|^2 \Delta f\dd v
\end{align*}
The last two terms cancel out, and we obtain \eqref{e:Fisher derivative heat equation}. For a more complete discussion of the Fisher information we recommend Villani's recent preprint \cite{villani2025}, which aside from a complete historical discussion also discusses the result of Imbert, Villani, and the second author for the Boltzmann equation \cite{ImbSilVil2024}.

While the computation we have just described is elementary, it seems to be very specific to the heat equation. Studying parabolic equations with variable coefficients, and establishing conditions on the coefficients that guarantee the monotonicity of the Fisher information is a much more difficult problem that is not well understood.

\subsection{Lifting the Fisher information}

Consider now the Fisher information functional over probability densities in $\mathbb{R}^d$ and $\mathbb{R}^{2d}$, namely
\begin{align*}
  i(f) & = \int_{\mathbb{R}^d}|\nabla_v \ln f(v)|^2f(v) \dd v\\
  I(F) & = \int_{\mathbb{R}^{2d}}|\nabla_{v,w}\ln F(v,w)|^2F(v,w)\ddvw
\end{align*}	
Since $\ln f\otimes g = \ln f(v) + \ln g(w)$ we have
\begin{align*}
  |\nabla \ln (f\otimes g)(v,w)|^2 = |\nabla \ln f(v)|^2+|\nabla \ln g(w)|^2
\end{align*}
If we multiply both sides by  $f(v)g(w)$ and integrate for $(v,w)$ in $\mathbb{R}^{2d}$ we have
\begin{align*}
  & \int_{\mathbb{R}^{2d}} |\nabla_{v,w} \ln(f(v)g(w))|^2f(v)g(w)\ddvw\\
  & = \int_{\mathbb{R}^{d}} |\nabla \ln f|^2f\dd v + \int_{\mathbb{R}^{d}} |\nabla \ln g|^2g\dd w
\end{align*}
This can be written more compactly as
\begin{align*}
  I(f\otimes g) = i(f)+i(g).
\end{align*}
This shows the Fisher information is additive with respect to tensor products of densities (befitting a quantity named information: two identical independent trials produce \emph{twice} the information).

Now consider an arbitrary probability density $F(v,w)$ in $\mathbb{R}^{2d}$ with marginals $f$ and $g$, namely 
\begin{align*}
  f(v) = \int_{\mathbb{R}^{d}}F(v,w)\dd w,\; g(w) = \int_{\mathbb{R}^{d}}F(v,w)\dd v
\end{align*}
From here, Jensen's inequality says that for any $v$
\begin{align*}
   |\nabla f(v)|^2 & = \left |\int_{\mathbb{R}^N}(\nabla_v\ln F(v,w))F(v,w)\dd w \right |^2\\
     & \leq f(v)\int_{\mathbb{R}^N}|\nabla_v \ln F(v,w)|^2F(v,w)\dd w
\end{align*}
and for any $w$
\begin{align*}
   |\nabla g(w)|^2 & \leq g(v)\int_{\mathbb{R}^N}|\nabla_w \ln F(v,w)|^2F(v,w)\dd v
\end{align*}
We thus have the pointwise inequalities
\begin{align*}
   f(v)|\nabla_v \ln f|^2 & \leq \int_{\mathbb{R}^N}|\nabla_v \ln F|^2F(v,w)\dd w\\
   g(w)|\nabla_w \ln g|^2 & \leq \int_{\mathbb{R}^N}|\nabla_w \ln F|^2F(v,w)\dd v   
\end{align*}
Integrating these two inequalities we arrive at 
\begin{align*}
  i(f)+i(g) \leq I(F)
\end{align*}
We state this as a lemma.

\begin{lemma} \label{l:subadditivity of Fisher}
Given a probability density $F$ in $\mathbb{R}^d\times \mathbb{R}^d$ with marginals $f,g$ in $\mathbb{R}^d$, we have
\begin{align*}
I(f\otimes g) \leq I(F)	
\end{align*}	
In particular, when $f=g$ we have 
\begin{align*}
2i(f) = I(f\otimes f) & \leq I(F)
\end{align*}	
\end{lemma}
While the proof of Lemma \ref{l:subadditivity of Fisher} is elementary, it has some remarkable consequences. This was first observed by Carlen \cite{Car1991}, who used it to obtain simple proofs of the Blachman-Stam inequality, and the characterization of the equality case in the logarithmic Sobolev inequality. In our context, Lemma \ref{l:subadditivity of Fisher} tells us that the Fisher information in $\mathbb{R}^{2d}$ is a natural lifted functional of the Fisher information in $\mathbb{R}^d$.

\section{The lifted equation}

We have established that the monotonicity of the Fisher information for the Landau equation boils down to the monotonicity of the Fisher information for the lifted equation. It is now time to focus our study on the lifted equation itself.  
\begin{align}\label{e:lifted equation}
  \partial_tF & = Q(F) \text{ in } (0,\infty)\times \mathbb{R}^6
\end{align}
An important observation about this equation is its solutions amount (in the right system of coordinates) to families of solutions to the heat equation on $\mathbb{S}^2$. 

\begin{figure}[h]

\vspace{-0.05in}
\centering
\begin{picture}(120,150)

\put(-30,00){\includegraphics[scale=0.25]{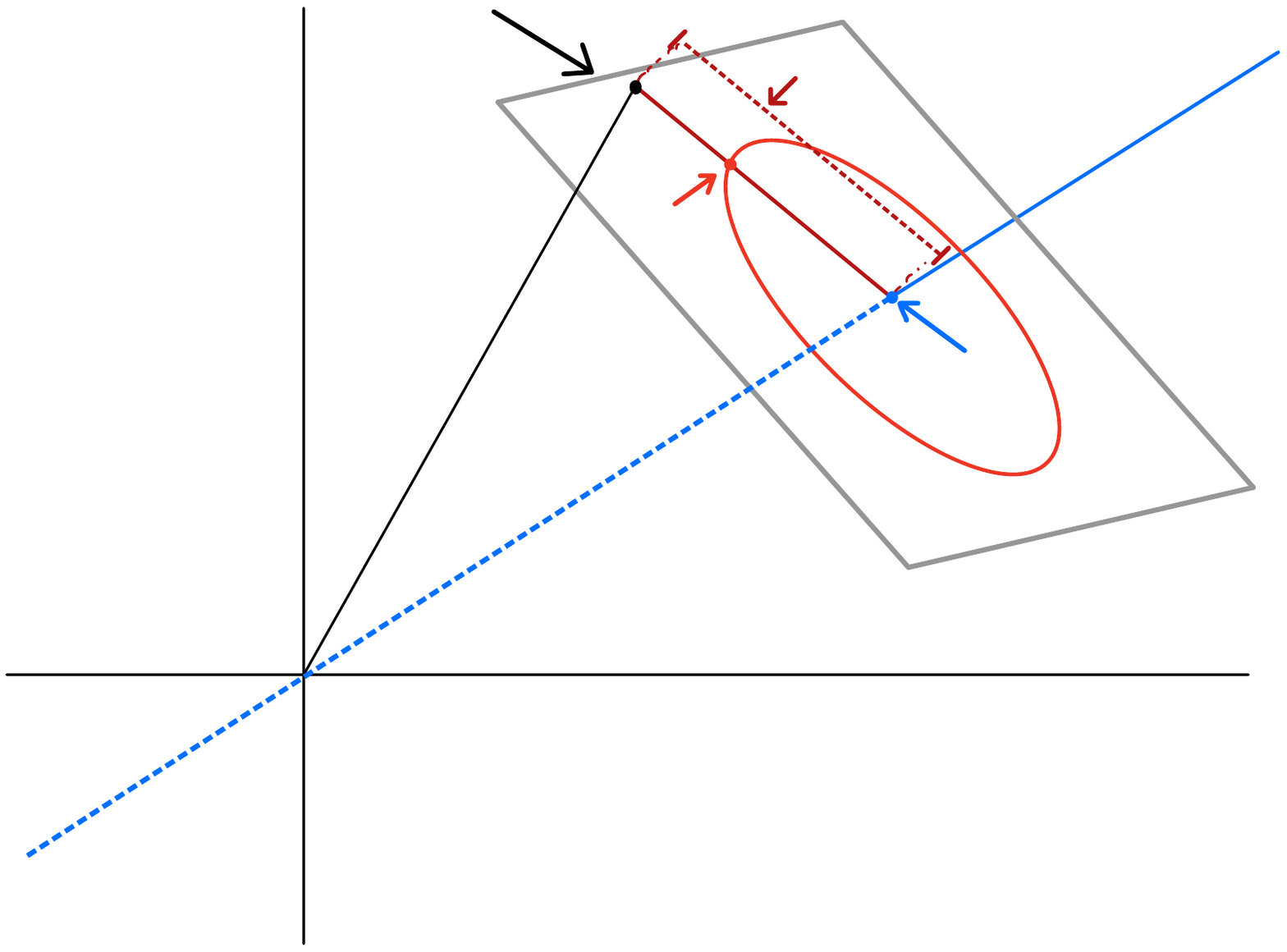}}
\put(40,10){$v = z + r\sigma$, $w = z - r\sigma$}
\put(105,86){{\color{blue}$(z,z)$}}
\put(144,120){{\color{blue}$\{v=w\}$}}
\put(18,145){$(v,w)$}
\put(65,110){{\color{red} $\sigma$}}
\put(92,136){{\color{red} $r$}}
\put(150,35){$v$}
\put(5,140){$w$}
\put(-45,140){$\mathbb{R}^3\times \mathbb{R}^3$}

\end{picture}

\vspace{-0.25in}
\caption{Change of variables from $(v,w)$ to $(z,r,\sigma)$.}\label{f:coordinates}
\end{figure}

To explain this we need to switch to a different system of coordinates (see Figure \ref{f:coordinates})
\begin{align*}
  (z,r,\sigma) \mapsto (z+r\sigma,z-r\sigma) \in \mathbb{R}^6	 
\end{align*}
Here $(z,r,\sigma) \in \mathbb{R}^3 \times (0,\infty) \times \mathbb{S}^2$. We shall write
\begin{align*}
  \bar F(t,z,r,\sigma) & := F(t,z+r\sigma,z-r\sigma)\\
  \bar F_{\text{in}}(z,r,\sigma) & := F_{\text{in}}(z+r\sigma,z-r\sigma)
\end{align*}
Going forward a bar will indicate a function of the $(z,r,\sigma)$ variables. The operator $Q$ takes a particularly simple form in these coordinates. Observe that the matrix $a_{ij}$, introduced in \eqref{e:aij}, is a scalar multiple of the projection onto the orthogonal complement of the vector $v-w$. In terms of the new coordinates, this means that $a_{ij}(v-w)$ is a scalar function times the projection onto the tangent space to the sphere $\mathbb{S}^2$ at $\sigma$. Effectively, it allows us to interpret the operator $Q$ as a scalar multiple of the Laplace-Beltrami operator on the sphere $\mathbb{S}^2$. The precise statement is given in the following lemma.

\begin{lemma}\label{l:lifted equation is heat equation on sphere}
  In terms of the coordinates $(z,r,\sigma)$, the operator $Q$ takes the form
  \begin{align*}
    \bar Q(F)(z,r,\sigma) = \alpha(r)\Delta_{\sigma}\bar F(z,r,\sigma)
  \end{align*}

  Consequently, the equation \eqref{e:lifted equation} is simply the heat equation with respect to $t$ and $\sigma$, with a diffusion coefficient $\alpha(r)$, for each fixed $(z,r)$. \begin{align}\label{e:lifted equation spherical variables}
    \partial_t \bar F = \alpha(r) \Delta_\sigma \bar F 
  \end{align}	  
\end{lemma}
We discussed in Section \ref{s:the Fisher information} how the Fisher information is monotone decreasing along the flow of the heat equation in Euclidean space. The computation we described in Section \ref{s:the Fisher information} is quite rigid. There is no easy criterium to determine when the Fisher information is monotone decreasing along the flow of a general linear parabolic equation. We have to use the very special form of $Q$ to analyze the monotonicity of the Fisher information for this particular diffusion. Lemma \ref{l:lifted equation is heat equation on sphere} suggests we should look separately at the Fisher information with respect to $z,r,$ and $\sigma$.  

\section{Fisher along layers}

The gradients $\nabla_{v,w}F$ and $\nabla_{z,r,\sigma}\bar F$ satisfy the relation
\begin{align*}
  |\nabla_{v,w}F|^2 = \tfrac{1}{2}\left (r^{-2}|\nabla_\sigma\bar F|^2 + (\partial_r \bar F)^2 + |\nabla_z \bar F|^2 \right )
\end{align*}
where the left hand side is evaluated at $(v,w) = (z+r\sigma,z-r\sigma)$ and the right hand side at $(z,r,\sigma)$. In light of this, we will decompose the Fisher information functional as the sum
\begin{align*}
  I(F) = \tfrac{1}{2}\left ( I_{\text{par}}(F) + I_{\text{sph}}(F)+I_{\text{rad}}(F)\right ),
\end{align*}
where 
\begin{align*}
  I_{\text{par}}(F) & := \int \frac{|\nabla_z\bar F(z,r,\sigma)|^2}{\bar F(z,r,\sigma)} r^{2} \ddsrz\\
  I_{\text{sph}}(F) & := \int \frac{|\nabla_\sigma \bar F(z,r,\sigma)|^2}{\bar F(z,r,\sigma)}\ddsrz\\
  I_{\text{rad}}(F) & := \int \frac{(\partial_r \bar F(z,r,\sigma))^2}{\bar F(z,r,\sigma)}r^{2}\ddsrz.
\end{align*}
These three functionals capture the components of the Fisher direction along the parallel directions ($z$), spherical directions ($\sigma$), and radial direction ($r$). We are going to estimate the derivative of each of these separately, resulting in the following inequalities
\begin{align*}
  \langle I_{\text{par}}'(F),Q(F)\rangle & \leq 0\\
  \langle I_{\text{sph}}'(F),Q(F)\rangle & = -2\int  \alpha\Gamma_{2}(\ln \bar F)\bar F \ddsrz\\
  \langle I_{\text{rad}}'(F),Q(F)\rangle & \leq \int ((\sqrt{\alpha})')^2  |\nabla_\sigma \ln \bar F|^2\bar F r^2\ddsrz 
\end{align*}
The operator $\Gamma_2$ denotes the \emph{iterated Carr\'e du champ operator} from Bakry-\'Emery theory \cite{bakry2006diffusions}, for the spherical Laplacian. It has a simple explicit form 
\begin{align}\label{e:Gamma deux}
  \Gamma_{2}(\ln \bar F) = |\nabla_\sigma^2\ln \bar F|^2+|\nabla_\sigma \ln \bar F|^2
\end{align}
Each inequality will follow in turn from an estimation of the invariance (or lack of) of the lifted equation \eqref{e:lifted equation spherical variables} with respect to shifts in the variables $z,\sigma,$ and $r$. Note in particular that from the above $I_{\text{par}}$ and $I_{\text{sph}}$ are monotone decreasing, but $I_{\text{rad}}$ may not if $\alpha' \neq 0$. The underlying reason for this is that the operator $Q$ is invariant under parallel shifts (translations in $z$) and rotations (rotations in $\sigma$). The operator $Q$ is not invariant under radial shifts unless $\alpha$ is constant (Maxwell molecules). The correspondence between symmetries of linear parabolic equations like the Fokker-Planck equation and monotonicity of Fisher information functionals has been known for some time in the physics literature, see for example the work of Plastino and Plastino \cite{Plastino1996}.

The lifting property of $i(f)$ and the estimates for the time derivatives of the three Fisher functionals combine into the following theorem.
\begin{thm}\label{t:Fisher time derivative estimate}
Writing $F = f\otimes f$, we have
\begin{align*}
\frac{d}{dt}i(f) & \leq  -\int \alpha \; \Gamma_{2}(\ln \bar F)\bar F\ddsrz \\
& \;\;\;\; + \int ((\sqrt{\alpha})')^2 |\nabla_\sigma \ln \bar F|^2\bar Fr^2\ddsrz
\end{align*}

\end{thm}

The first term on the right hand side is non-positive, and the second non-negative. In order to conclude that the Fisher information is monotone decreasing along the flow of the Landau equation, we must show the second term controls the first. We will show this for each value of $z$ and $r$, as a consequence of an inequality for functions on the sphere, in the family of the log-Sobolev inequality.

Note that when $\alpha$ is constant the second term vanishes and Theorem \ref{t:Fisher time derivative estimate} immediately implies the monotonicity of the Fisher information. This simpler scenario corresponds to the Maxwell molecules case $\alpha \equiv 1$, which is the case that had been previously studied in \cites{Mck1966,Tos1992,Vil1998}

In the next three subsections, we justify the formulas for the derivatives of $I_{\text{par}}$, $I_{\text{sph}}$, and $I_{\text{rad}}$, leading to Theorem \ref{t:Fisher time derivative estimate}. To explain these computations, we use a characterization of the Fisher information in terms of the Kullback-Leibler divergence (KL), also called \emph{relative entropy}. The relative entropy measures the difference between two probability distributions by
\begin{align*}
\textnormal{KL}(\bar F,\bar G) := \int_{\mathbb{S}^2}\bar F\ln(\bar F/\bar G)\dd\sigma
\end{align*}
Then, the Fisher information itself measures the change of a density with respect to this metric, 
\begin{align}\label{e:from KL to Fisher}
  {\frac{d^2}{d h^2}}_{\mid_{h=0}}\textnormal{KL}(F_{\theta},F_{\theta+h}) = \frac{1}{2}\int_{\mathbb{S}^2} \frac{(\partial_\theta F_\theta)^2}{F_\theta}\dd \sigma
\end{align}
The heat flow on the sphere $\mathbb{S}^2$ is contractive with respect to the relative entropy. Accordingly, if $\bar F,\bar G$ each solve \eqref{e:lifted equation spherical variables}, then for every $z,r$ fixed we have
\begin{align}\label{e:derivative of relative entropy}
\frac{d}{dt}\textnormal{KL}_{z,r}(\bar F,\bar G) = -\alpha \int_{\mathbb{S}^2} \bar F|\nabla_\sigma \ln (\bar F/\bar G)|^2\dd\sigma
\end{align}

We now proceed to compute the time derivative of the Fisher information functionals $I_{\text{par}}$, $I_{\text{sph}}$, and $I_{\text{rad}}$ by using \eqref{e:from KL to Fisher} in combination with \eqref{e:derivative of relative entropy}.

\subsection{Fisher along parallel directions}
\;
\begin{figure}[h]

\vspace{-0.32in}
\begin{center}
\includegraphics[scale=0.23]{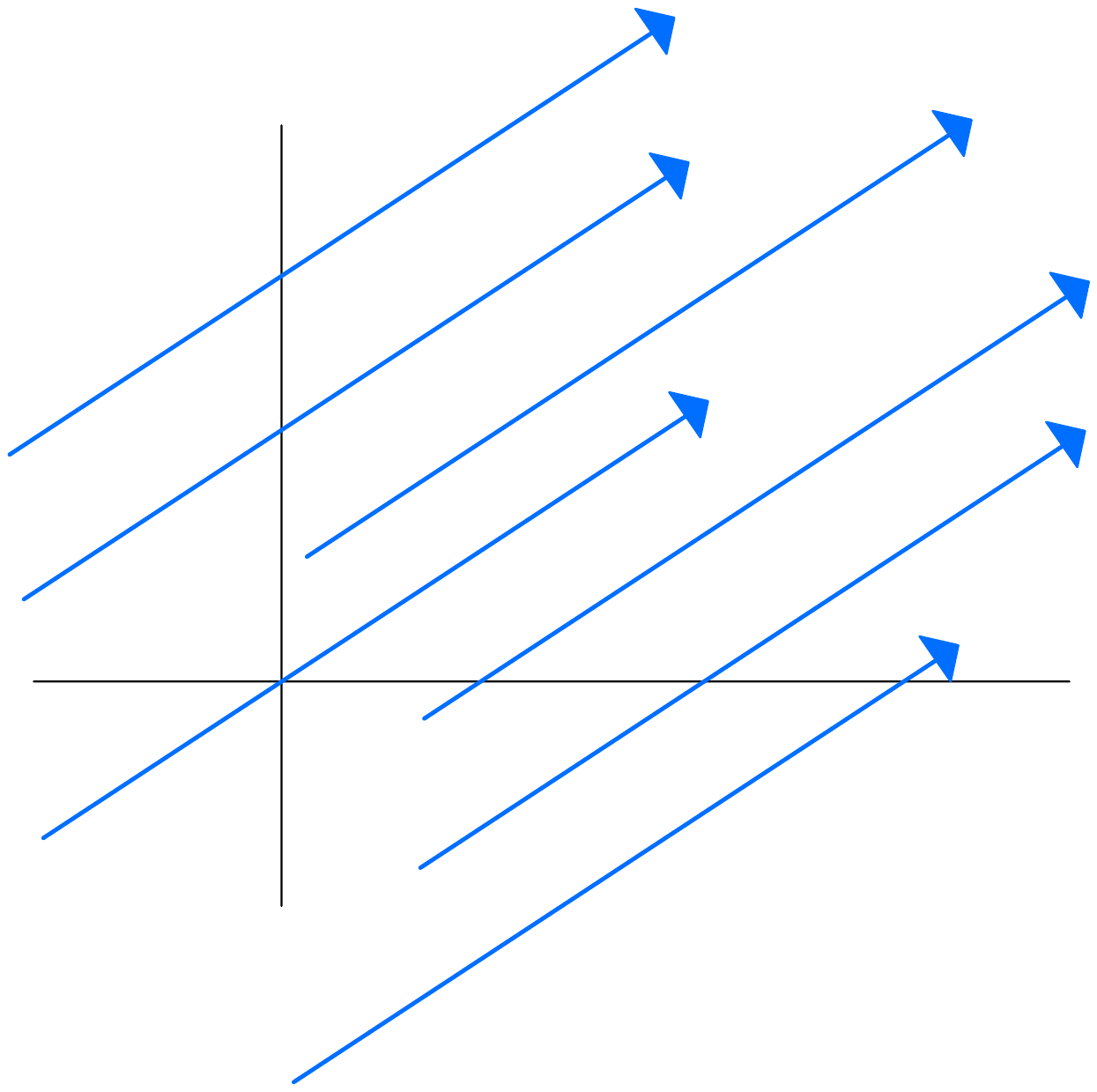}

\end{center}

\vspace{-0.32in}

\end{figure}

We start with the estimate for parallel directions. 
\begin{align*} 
  \langle I_{\text{par}}'(F),Q(F)\rangle \leq 0.
\end{align*}
For $e\in \mathbb{R}^3$ and $h>0$ we consider the function
\begin{align*}
  \bar G_h(t,z,r,\sigma) = \bar F(t,z+he,r,\sigma) 
\end{align*}
Note that since $\bar G_0 = \bar F$, $\textnormal{KL}_{z,r}(\bar F,\bar G_h)$ attains its minimum equal to $0$ at $h=0$. As we mentioned in \eqref{e:from KL to Fisher}, the second derivative of $\textnormal{KL}_{z,r}(\bar F,\bar G_h)$ with respect to $h$ at $h=0$ is equal to the Fisher information in the direction of $e$. Thus, for each $z,r$, we have
\begin{align*}
\lim \limits_{h\to 0^+}\frac{1}{h^2}\textnormal{KL}_{z,r}(\bar F,\bar G_h) = \frac{1}{2}\int_{\mathbb{S}^2} \frac{(\nabla_z \bar F\cdot e)^2}{\bar F}\dd\sigma
\end{align*}

To show that the left hand side is monotone decreasing with respect to $t$, we show that for every $h \neq 0$, $\textnormal{KL}_{z,r}(\bar F,\bar G_h)$ is a non-increasing function of $t$. This stems from the fact that both $\bar F$ and $\bar G_h$ solve the lifted equation \eqref{e:lifted equation} and so \eqref{e:derivative of relative entropy} applies. Dividing the identity by $h^2$ and taking $h\to 0^+$ it follows that
\begin{align*}
& \frac{1}{2}\frac{d}{dt}\int_{\mathbb{S}^2} \frac{(\nabla_z \bar F\cdot e)^2}{\bar F}\dd\sigma \\
& = -\alpha \lim\limits_{h\to 0^+}\int_{\mathbb{S}^2} \bar F|\tfrac{1}{h}\nabla_\sigma \ln (\bar F/\bar G_h)|^2\dd\sigma \leq 0
\end{align*}
It is straightforward to show the limit on the right exists but all we will need is that it has the right sign. Combining the resulting inequality for $e$'s in an orthonormal basis we conclude that for each $z,r$
\begin{align*}
& \frac{1}{2}\frac{d}{dt}\int_{\mathbb{S}^2} \frac{|\nabla_z \bar F|^2}{\bar F}\dd\sigma \leq 0.
\end{align*}
Integrating in $r^2\mathrm{d}z\mathrm{d}r$ the estimate for $I_{\text{par}}$ follows.

\subsection{Fisher along spherical directions}

\begin{figure}[h]
	
\vspace{-0.225in}
\begin{center}
\includegraphics[scale=0.23]{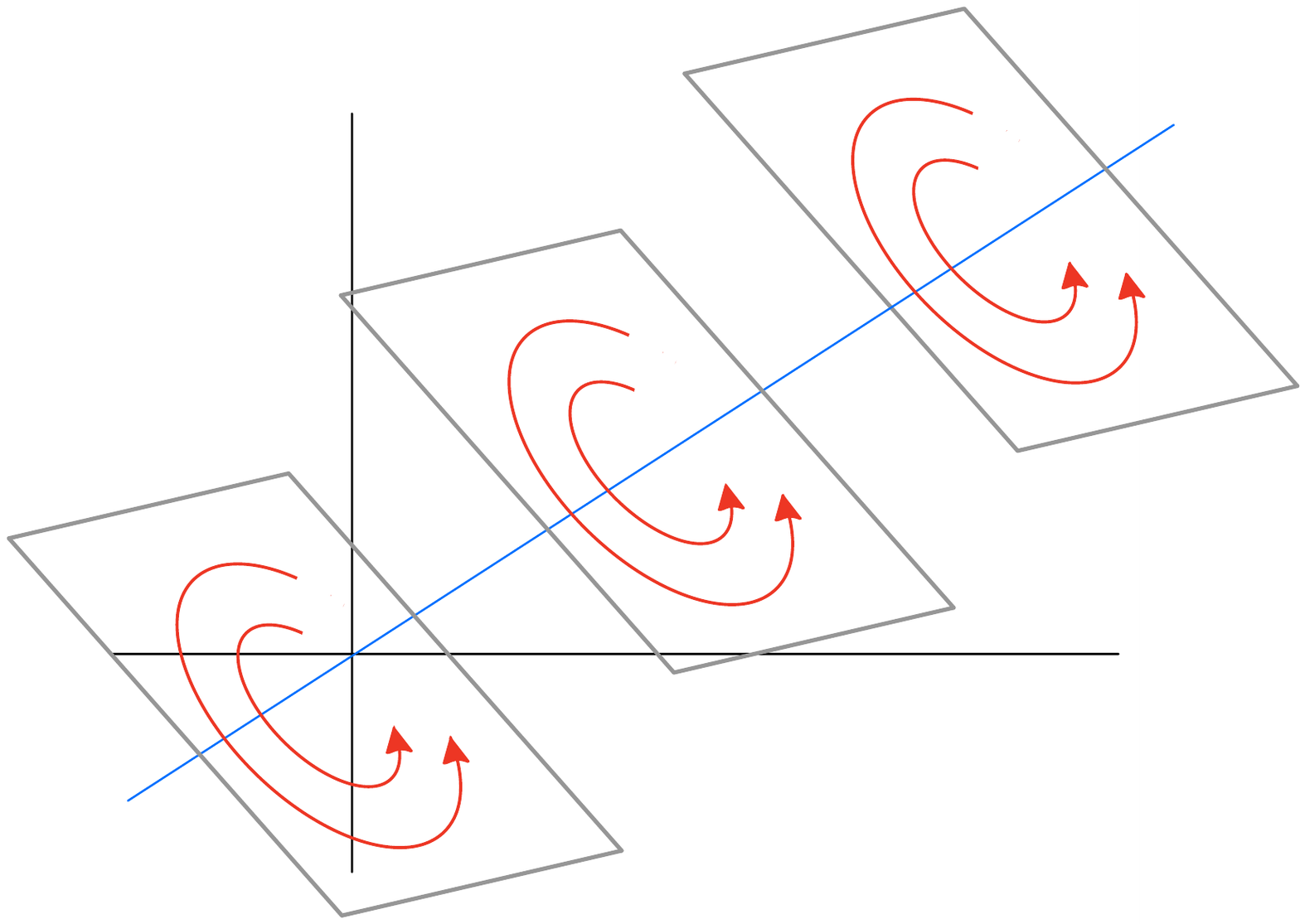}

\end{center}

\vspace{-0.32in}

\end{figure}

Now we analyze the spherical directions to show
\begin{align*}
  \langle I_{\text{sph}}'(F),Q(F)\rangle = -2\int \alpha \Gamma_{2}(\ln \bar F)\bar F \ddsrz,
\end{align*}
where $\Gamma_2$ is as in \eqref{e:Gamma deux}. Take a one-parameter family of rotations $h\mapsto R_h$. The function
\begin{align*}
  \bar G_h(t,z,r,\sigma) = \bar F(t,z,r,R_h\sigma) 
\end{align*}
solves \eqref{e:lifted equation spherical variables} thanks to rotational invariance. Denoting by $b$ the vector field generating $R_h$, \eqref{e:from KL to Fisher} says that
\begin{align*}
\lim \limits_{h\to 0^+}\frac{1}{h^2}\textnormal{KL}_{z,r}(\bar F,\bar G_h) = \frac{1}{2}\int_{\mathbb{S}^2} \frac{(\nabla_\sigma \bar F\cdot b)^2}{\bar F}\dd\sigma
\end{align*}
Then, dividing by $h^2$ and taking $h\to 0^+$ we have
\begin{align*}
& \frac{1}{2}\frac{d}{dt}\int_{\mathbb{S}^2} \frac{(\nabla_\sigma \bar F\cdot b)^2}{\bar F}\dd\sigma \\
& = -\alpha \lim\limits_{h\to 0^+}\int_{\mathbb{S}^2} \bar F|\tfrac{1}{h}\nabla_\sigma \ln (\bar F/\bar G_h)|^2\dd\sigma
\end{align*}
A bit of basic calculus on surfaces shows that
\begin{align*}
\lim\limits_{h\to 0}\tfrac{1}{h}\nabla_\sigma \ln (\bar F/\bar G_h) = (\nabla_\sigma^2 \ln \bar F )b + (\nabla_\sigma \ln \bar F\cdot b)\sigma 
\end{align*}
with the Hessian and gradient on $\mathbb{S}^2$. We arrive at
\begin{align*}
& \frac{1}{2}\frac{d}{dt}\int_{\mathbb{S}^2} \frac{(\nabla_\sigma \bar F\cdot b)^2}{\bar F}\dd\sigma \\
& = -\alpha \int_{\mathbb{S}^2} \bar F(|(\nabla_\sigma^2 \ln \bar F )b|^2+ |(\nabla_\sigma \ln \bar F\cdot b)|^2)\dd\sigma
\end{align*}
We average this over all vector fields $b$ corresponding to generators of rotations and arrive at
\begin{align*}
& \frac{d}{dt}\int_{\mathbb{S}^2} \frac{|\nabla_\sigma \bar F|^2}{\bar F}\dd\sigma = -2\alpha \int_{\mathbb{S}^2} \Gamma_2(\ln \bar F)\bar F\dd\sigma
\end{align*}
Integrating in $z,r$ we obtain the desired identity. 

\subsection{Fisher along the radial direction}

\begin{figure}[h]

\vspace{-0.225in}

\begin{center}
\includegraphics[scale=0.23]{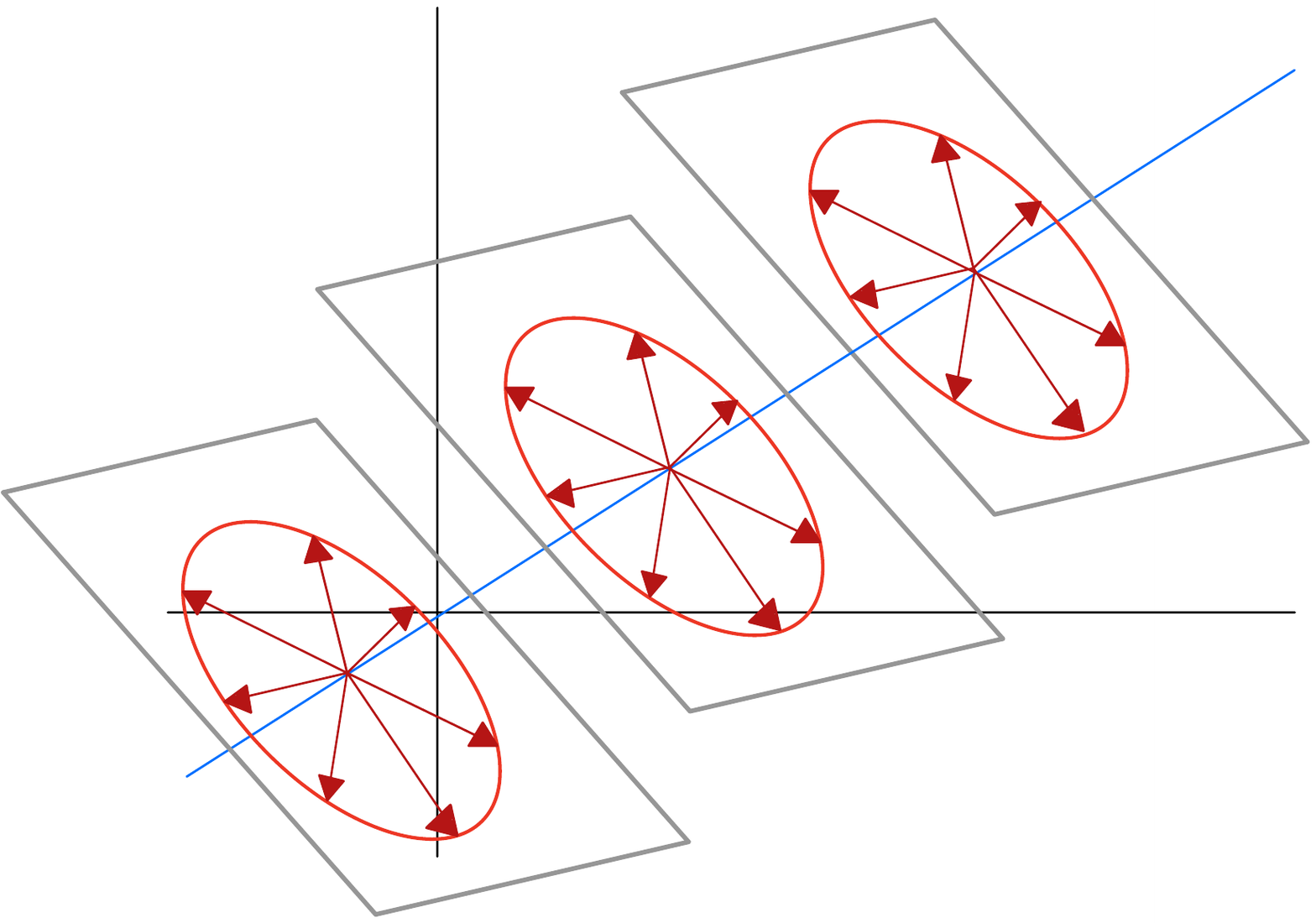}
\end{center}

\vspace{-0.32in}

\end{figure}

All that remains is the radial direction estimate
\begin{align*}
  \langle I_{\text{rad}}'(F),Q(F)\rangle & \leq 2\int ((\sqrt{\alpha})')^2|\nabla_\sigma \ln \bar F|^2\bar F r^2\ddsrz
\end{align*}
This estimate captures the one problematic term in Theorem \ref{t:Fisher time derivative estimate}, as is the one that could potentially drive $I(F)$ towards increasing.  For this term we use
\begin{align*}
\bar G_h(t,z,r,\sigma) = \bar G_h(t,z,r+h,\sigma)
\end{align*}
and, as the reader may now expect we have
\begin{align*}
  \lim \limits_{h\to 0}\frac{1}{h^2}\textnormal{KL}_{z,r}(\bar F,\bar G_h) = \frac{1}{2}\int_{\mathbb{S}^2} \frac{(\partial_r \bar F)^2}{\bar F}\dd\sigma
\end{align*}
Given that in general $\alpha(r) \neq \alpha(r+h)$ we cannot say anymore that $\bar F$ and $\bar G_h$ solve the same equation, so \eqref{e:derivative of relative entropy} no longer holds. Now we must somehow measure the extent to which the equation for $\bar F$ is no longer invariant under shifts in the parameter $r$. For this reason we write the equation for $\bar G_h$ as follows
\begin{align*}
\partial_t \bar G_h = \alpha(r) \Delta_\sigma \bar G_h + (\delta_h \alpha) \Delta_\sigma \bar G_h
\end{align*}
where $\delta_h \alpha := \alpha(r+h)-\alpha(r)$. Taking into account this correction term we have
\begin{align*}
\frac{d}{dt}\textnormal{KL}_{z,r}(\bar F,\bar G_h) & = -\alpha \int_{\mathbb{S}^2} \bar F|\nabla_\sigma \ln (\bar F/\bar G_h)|^2\dd\sigma\\
  & \;\;\;\; - (\delta_h \alpha)\int_{\mathbb{S}^2} (F/G_h)\Delta_\sigma G_h\dd\sigma
\end{align*}
As before we divide by $h^2$ and take $h\to0^+$, yielding
\begin{align*}
\frac{1}{2}\frac{d}{dt}\int_{\mathbb{S}^2} \frac{(\partial_r \bar F)^2}{\bar F}\dd\sigma & = -\alpha\int_{\mathbb{S}^2}|\nabla_\sigma \partial_r\ln\bar F|^2\bar F\dd\sigma\\
  & \;\;\;\;+\alpha'\int_{\mathbb{S}^2}\partial_r\ln \bar F\Delta_\sigma \bar F \dd\sigma
\end{align*}
Integrating by parts in the second term we obtain 
\begin{align*}
  & -\int_{\mathbb{S}^2}\alpha|\nabla_\sigma\partial_r\ln \bar F|^2 \bar F +\alpha' (\nabla_\sigma \partial_r \ln \bar F)\cdot \nabla_{\sigma}\bar F\dd\sigma
\end{align*}
The integrand involves an incomplete square, completing this square we rewrite the integrand as
\begin{align*}
-|\partial_r(\sqrt{\alpha}\nabla_\sigma\ln \bar F)|^2\bar F +((\sqrt{\alpha})')^2|\nabla_\sigma \ln \bar F|^2\bar F 
\end{align*}
Dropping the first (non-positive) term we arrive at
\begin{align*}
  \frac{1}{2}\frac{d}{dt}\int_{\mathbb{S}^2} \frac{(\partial_r \bar F)^2}{\bar F}\dd\sigma \leq \int_{\mathbb{S}^2}((\sqrt{\alpha})')^2|\nabla_\sigma \ln \bar F|^2 \dd\sigma
\end{align*}
Integration respect to $r^2\mathrm{d}z\mathrm{d}r$ yields the inequality.

\section{An inequality on $\mathbb{S}^2$}

From Theorem \ref{t:Fisher time derivative estimate} we see $\frac{d}{dt}i(f)$ will be $\leq 0$ provided that for any fixed $z,r$ the following inequality holds 
\begin{align*}
\frac{(r\alpha'(r))^2}{4\alpha(r)}\int_{\mathbb{S}^2}\frac{|\nabla_\sigma \bar F|^2}{\bar F}\dd\sigma \leq \alpha \int_{\mathbb{S}^2} \Gamma_{2}(\ln \bar F)\bar F\dd\sigma
\end{align*}
We recall that $F(t,v,w) = F(t,w,v)$ and observe that in terms of the $(z,r,\sigma)$ variables this symmetry means $\bar F$ is even in $\sigma$, namely $\bar F(z,r,\sigma) = \bar F(z,r,-\sigma)$. 

Therefore, our goal will follow if we show that for any smooth positive function $\bar f:\mathbb{S}^2\to\mathbb{R}$ satisfying the symmetry $\bar f(\sigma)=\bar f(-\sigma)$ we have the inequality
\begin{align}\label{e:BakryEmery condition}
\lambda \int_{\mathbb{S}^2}|\nabla_\sigma \ln (\bar f)|^2 \bar f\dd\sigma \leq \int_{\mathbb{S}^2} \Gamma_{2}(\ln \bar f)\bar f\dd\sigma
\end{align}
for $\lambda>0$ independent of $\bar f$ such that $\lambda \geq \Lambda(\alpha)$, where
\begin{align*}
  \Lambda(\alpha) := \left (\sup_{r>0}\frac{r\alpha'(r)}{2\alpha(r)}\right )^2.
\end{align*}
Since $\Gamma_2(\ln \bar f)= |\nabla_\sigma^2 \ln\bar f|^2 + |\nabla_\sigma \ln \bar f|^2$ we have
\begin{align*}
 |\nabla_\sigma \ln \bar f|^2\bar f \leq \Gamma_2(\ln \bar f)\bar f \text{ in } \mathbb{S}^2
\end{align*} 
and so it follows trivially that \eqref{e:BakryEmery condition} holds with $\lambda=1$. Already from here one obtains the monotonicity of the Fisher information for $\alpha(r) = r^\gamma$ with $\gamma \in [-2,2]$. This does not cover the Coulomb case $\gamma=-3$. A more detailed analysis leads to the following lemma.
\begin{lemma}\label{l:inequality on the sphere}
  There is a $\lambda\geq \frac{11}{2}$ such that 
  \begin{align*}	
  \lambda \int_{\mathbb{S}^2}|\nabla_\sigma \ln f|^2f \dd\sigma \leq \int_{\mathbb{S}^2}\Gamma_{2}(\ln f)f \dd\sigma
  \end{align*}
  for any $f:\mathbb{S}^2\to (0,\infty)$, smooth and even.
\end{lemma}

In \cite{GuiSil2023} this inequality was proved with a weaker estimate on the best constant $\lambda$, namely $\lambda \geq \tfrac{19}{4}$. The improved bound $\lambda \geq \tfrac{11}{2}$ was obtained by Ji in \cite{Ji2024}, extending the range of potentials for which the monotonicity of $i(f)$ holds. It is also shown in \cite{Ji2024} that the $\lambda$ cannot be larger than $\approx 5.73892$.  

The inequality \eqref{e:BakryEmery condition} is known as the Bakry-\'Emery $\Gamma_2$ criterion \cite{bakry2006diffusions} for log-Sobolev inequalities. The best constant $\lambda$ for $S^2$, without the symmetry assumption on $\bar f$, is well-known to be $\lambda=2$. The observation that our functions of interest are symmetric allows us to improve the constant $\lambda$ in \eqref{e:BakryEmery condition} to $\tfrac{11}{2}$, allowing us to treat the Coulomb case $\gamma=-3$.

Thanks to Theorem \ref{t:Fisher time derivative estimate} and Lemma \ref{l:inequality on the sphere} we conclude $i(f(t))$ is always decreasing in time as long as the potential $\alpha(r)$ satisfies the bound $\Lambda(\alpha) \leq \tfrac{11}{2}$. For $\alpha(r)=r^{\gamma}$, this condition becomes
\begin{align*}
  |\gamma| \leq \sqrt{22}
\end{align*}
which covers the case $\gamma=-3$ of Coulomb potentials.

\section{Conclusion}

We explained the ideas involved in the proof of the monotonicity of the Fisher information for the space-homogeneous Landau equation. It is not a frequent event to find new monotone quantities for well-known equations in mathematical physics. A control of a higher-order quantity naturally says a lot about the behavior of the equation. In particular, it allows us to prove that a smooth solution exists for all time in the case of the homogeneous Landau \cite{GuiSil2023}. 

The ideas described in this article are applied to the Boltzmann equation in \cite{ImbSilVil2024}. In that case, the lifted operator $Q$ is an integro-differential diffusion operator with respect to the spherical variable $\sigma$. The analysis is consequently more involved. There is a major new difficulty for the case of the Boltzmann equation, which is to develop a new Bakry-\'Emery criterion, with explicit constants, for integro-differential diffusion operators on the sphere.

For now, the techniques described in this article are not applicable to inhomogeneous kinetic equations. Some work remains to be done in this direction.

We also hope that this method will promote further study of the best constants in the Bakry-\'Emery criterion for log-Sobolev inequalities on the sphere, with more general diffusion operators including integro-differential ones.

\bibliographystyle{plain}
\bibliography{landau_notices}

\end{document}